\documentclass[a4paper,reqno]{amsart}
\usepackage{latexsym}
\usepackage{amssymb,amsthm,amsmath}
\usepackage[dvips]{graphicx}
\usepackage{xcolor}
\usepackage{array}
\usepackage{hyperref} 
\usepackage{stackengine}
\usepackage{amsaddr}

\usepackage{graphicx}
 \usepackage{subcaption}
 \usepackage{epsfig}
 \usepackage[latin1]{inputenc}
 \usepackage{a4wide}
 \usepackage{caption}
 \usepackage{comment}
\usepackage{tcolorbox}   
\definecolor{mycolor}{rgb}{0.122, 0.435, 0.698}
\newcommand{\mybox}[1]{%
  \begin{tcolorbox}[colframe=mycolor,boxrule=0.2pt,arc=4pt,left=0pt,right=0pt,top=0pt,bottom=0pt,boxsep=2pt] #1
  \end{tcolorbox}%
}

\calclayout

\theoremstyle{plain}
\newtheorem{theorem}{Theorem}

\newtheorem{example}[theorem]{Example}

\theoremstyle{definition}

\theoremstyle{remark}

\reversemarginpar

\newcounter{zd}

\newcounter{zdr}[subsection]

\def\pa{\partial}

\def\cal{\mathcal}

\def\R{\mathbb{R}}    

\begin{document}

\title[Fair coalition agreement]{Pre-electoral coalition agreement from the Black-Scholes point of view}

\author{D.~Mitrovi\'{c}}
\address{
University of Vienna, Faculty of Mathematics, Oskar Morgenstern Platz-1,
1090 Vienna, Austria and}
\address{
University of Montenegro, Faculty of Mathematics and Natural Sciences, Cetinjski put bb,
81000 Podgorica, Montenegro}
\address{e-mail: darko.mitrovic@univie.ac.at}



\begin{abstract}
A political party can be considered as a company whose value depends on the voters support i.e. on the percentage of population supporting the party. Dynamics of the support is thus as a stochastic process with a deterministic growth rate perturbed by a white noise modeled through the Wiener process. This is in an analogy with the option modeling where the stock price behaves similarly as the voters' support. While in the option theory we have the question of fair price of an option, the question that we ask here is what is a reasonable level of support that the coalition of a "major" party (safely above the election threshold) and a "minor" party (under or around the election threshold) should achieve in order for the "minor" party to get one more representative.  
We shall elaborate some of the conclusions in the case of recent elections in Montenegro (June, 2023) which are particularly interesting due to lots of political subjects entering the race.
\end{abstract}

\maketitle

\section{Introduction}\label{sec:intro}

A defining trait of the majority of parliamentary democracies is that no single party can independently establish a government; instead, they are compelled to enter into some form of coalition. These coalitions can be formed either before or after the elections and from the perspective of political theory, there exists a wide array of research on this subject. Notably, post-electoral coalitions have undergone more extensive theoretical investigation. Conversely, there has been a noticeable dearth of research concerning pre-electoral coalitions, as noted by Golder \cite{golder}. In this present contribution, our aim is to introduce a mathematical approach to understanding the phenomenon of pre-electoral coalition formation.

From the purely mathematical point of view i.e. given the premise that voters will support a party regardless of whether it enters a coalition or not, it becomes evident that opting for a pre-electoral arrangement is a more logical choice. This is because, during the process of seat distribution, a political party inevitably loses a portion of its support (e.g. due to error in rounding -- see the Arrows paradox \cite{Arrows} for a comprehensive explanation). In other words, it would be more lucrative for  political parties with similar programs to enter a pre-electoral coalition since the subjects with higher percentage of votes are not only in absolute, but also in relative advantage with respect to those with lower percentage. The example from Table 1 from recent elections in Montenegro serves as an illustration of the latter situation. The resources throughout the paper are pre-electoral polls by the Center for Democracy and Human Rights (CEDEM) \cite{CEDEM} and the State Elections Committee of Montenegro \cite{DIK}. 

\begin{table}

\begin{tabular}{ c  c  c}

\begin{tabular}{ | m{5em} | m{2cm}| m{2cm}| } 
  \hline
   winner parties and JfA & actual mandates & percentage \\ 
  \hline
  MEN & 24 & 25.55\% \\ 
  \hline
  DPS & 21 & 23.26\% \\ 
  \hline 
  FfM & 13 & 14.76\% \\
  \hline
  JfA & 0 & 2.77\% \\
\hline
\end{tabular} & 
\begin{tabular}{ | m{5em} | m{3cm}|} 
  \hline
   winner parties & potential coalition   \\ 
  \hline
  MEN & 23  \\ 
  \hline
  DPS & 20  \\ 
  \hline 
  FfM+JfA & 15 \\
\hline
\end{tabular}

\end{tabular}
\caption{Left panel are actual results of the three "major" subjects Movement Europe Now (MEN), coalition around Democratic Party of Socialist (DPS) and coalition around For the Future of Montenegro (FfM). Results of the "minor" party Justice for All (JfA) are also provided. On the right panel are the results of the same entities if FfM would be in coalition with the "minor" party Justice For All (JfA).}
\end{table}


However, in the case when a party safely overcomes an election threshold or when a party does not have fairly compatible political partners, then it might be reasonable not to make any pre-electoral coalition. To this end, we shall call such parties the "major" parties while those parties that are below or around the election threshold we call the "minor" parties.

More specifically, when it comes to two major political parties with a secure support base above the election threshold, it becomes challenging to distribute seats in a fair manner. Unless there is a critical situation in the country, such as a long-term leader moving towards a dictatorship or, more generally, an extremely polarized society (see \cite{golder} for an in depth analysis), it seems an excessive risk for both of them to form a coalition. There are a couple of reasons for this hesitation. Firstly, some voters may be unwilling to vote for a coalition as they strongly support only one party, even though the ultimate outcome would be the same - a coalition government formed after the elections. Secondly, the more complex the coalition structure, the higher the likelihood of potential scandals arising, which would diminish the popularity of both parties. The example provided in Table 2 with data taken from the last two Montenegrin elections demonstrates such a situation.


\vspace{0.5cm}
\begin{table}

\begin{tabular}{ c  c }

\begin{tabular}{ | m{3cm} | m{2cm}| } 
  \hline
   no coalition (2020) & mandates  \\ 
  \hline
  URA & 4  \\ 
  \hline
  DEM & 10  \\ 
  \hline
\end{tabular} &
\begin{tabular}{ | m{3cm} | m{2cm}| } 
  \hline
   in coalition (2023) & mandates   \\ 
  \hline
  URA & 4  \\ 
  \hline
  DEM & 7  \\ 
  \hline
\end{tabular}

\end{tabular}
\caption{Left panel are results of the parties United Reform Action (URA) and Democrats (DEM) while acting independently in 2020. On the right panel are the results of the same subjects being in coalition in 2023.} 
\end{table} 


We are particularly interested in studying the role of "minor" parties, which refers to those whose support is around or below the election threshold. It is evident that these parties have a vested interest in securing a parliamentary status through pre-electoral coalitions. Simultaneously, it is also advantageous for a "major" party with similar political preferences to form a coalition to avoid potential loss of votes. The objective of our research paper is to propose a model that can assist parties in determining a fair pre-coalition agreement. Specifically, we aim to determine what percentage of support (considered and modeled as the price of an option) a "major" party should transfer to their "minor" coalition partner, thereby enabling the "minor" party to gain another representative in the parliament.

To provide a clearer explanation, let's examine the results of the last Montenegrin elections in June 2023, as presented in Table 1. Since the election threshold in Montenegro is 3$\%$ (for two representatives), we see that if the coalition "For the future of Montenegro" (FfM) were to transfer 0.23$\%$ of their support to the coalition "Justice for all," the latter would qualify for the parliament with two representatives. This transfer would not cost "For the future of Montenegro" a representative and would bring them two representatives of the closely aligned ally. In this scenario, a "transfer" of 0.23$\%$ is both beneficial and fair. We observe, nonetheless, that the potential representation of the "minor" party within a coalition is contingent upon pre-electoral polls. To illustrate, if 1.5$\%$ constitutes the threshold for earning one representative and the pre-electoral support for the "minor" party stands at 2$\%$, the coalition agreement must ensure at least one representative for the "minor" party. The question at hand is to determine the specific conditions under which the "minor" party would be entitled to an additional mandate i.e.:

\mybox{What level of support should a coalition of a "major" party and a "minor" party achieve in order for the minor party to gain an additional mandate?}

%

To this end, our starting point is the theory of options pricing, which aims to estimate the value of a financial contract in a just and fair manner. Common examples of such contracts are European put or call options, which provide the holder with the right, but not the obligation, to sell or buy the underlying financial instrument, respectively. It is important to note that the freedom to exercise or not exercise the agreement from the option contract has a cost and for many years, extensive research has been conducted to determine an appropriate price for entering into such contracts. This research culminated in the discovery of the Black-Scholes equation \cite{BS}, for which F. Black and M. Scholes were awarded the Nobel Prize in Economics. Given that we are dealing with at least two parties (two parties planning to form a coalition), we will require a variant of recently introduced multi-dimensional variant of the Black-Scholes equation \cite{Gui} to accurately assess the situation. We address a reader to \cite{Gui} for more information on the Black-Scholes equation.  

What we assume is that support to a political party behaves in the same way as the price of a company stock. To be more precise, let us consider a party say $P$ and denote by $S$ percentage of population supporting the party. We expect that the support grow/decrease by a steady rate $\mu$ while, using historical data, we can determine a fluctuation of the support (it is called volatility in finances) which we denote by $\sigma$. The situation here is somewhat simpler than when considering the option pricing since in the latter case, one is interested in the relative change off the stock price $S$ which leads to the following SDE
$$
\frac{dS}{S}=\mu dt+\sigma dW_t
$$ where $W_t$ is the Wiener process. Here, we are interested in the absolute support to the party and the starting point is the SDE
\begin{equation}
\label{1}
dS=\mu dt+\sigma dW_t.
\end{equation}  In the context of maximizing financial gains, the only relevant factor is the percentage change in price. In this scenario, investors naturally gravitate towards companies with the highest growth rates, as they yield the largest profits, regardless of the actual stock price. However, in the realm of politics, the absolute percentage of support is the key factor of interest. This is because it determines the number of representatives a party will have, irrespective of the relative growth in support compared to other parties.

Building upon equation \eqref{1}, we aim to derive an equation that provides information on the total number of votes a coalition should secure in order for a "minor" party to gain an additional mandate. 

%

The paper is organized as follows. 

After the introduction, we will provide a brief overview of the mathematical background and develop a model that pertains to the "price" of a pre-electoral coalition. In the final section, we will analyze the model using the recent elections in Montenegro as a case study. It is important to note that this analysis does not delve into the political landscape or parties in Montenegro, but rather serves as a theoretical exploration of the potential implications of applying the model.


\section{Mathematical background and the pre-elections coalition model}

We consider a fixed complete probability space $(\Omega, \{{\cal F}_t\}, {\cal F}, {\bf P})$, $t\in [0,T]$, $T>0$, with the sample space $\Omega$, the $\sigma$-algebra ${\cal F}$, the natural filtration $\{{\cal F}_t\}_{t\in [0,T]}$ generated by the Wiener process $W_t$, and the probability measure ${\bf P}$.

We also recall that the Wiener process $W_t$ is a real-valued stochastic process such that $W_0=0$, $W_t$ has independent increments, it has Gaussian increments $W_{t+u}-W_t \sim {\cal N}(u,0)$, and it has almost surely continuous paths. By $dW_t$ we denote the Wiener measure which is here simply an infinitesimal element of the It\^o integral $\int_0^t X_s dW_s:=\lim\limits_{n\to \infty} \sum\limits_{k=0}^nX_{t_k}(W_{t_{k+1}}-W_{t_{k}})$ for an exhausting partition of the set $[0,t]$ and $X_s$ being adapted with respect to the filtration $\{{\cal F}_s\}_{s\in [0,T]}$.

Let us finally recall the It\^o formula. Assume that the stochastic processes $S_i$, $i=1,2$, satisfy equations of type \eqref{1} with the constants $\mu=\mu_i$ and $\sigma=\sigma_i$, $i=1,2$. For the composition of a function $V\in C^1([0,T]\times \R^2)$ and stochastic processes $S_1$ and $S_2$, we have
\begin{equation}
\label{ito}
\begin{split}
dV(t,S_1,S_2)&=\left( \frac{\pa V}{\pa t} + \mu_1 \frac{\pa V}{\pa S_1}+\mu_2 \frac{\pa V}{\pa S_2}+\frac{\sigma_1^2}{2} \frac{\pa^2 V}{\pa S_1^2}+{\sigma_1\sigma_2} \frac{\pa^2 V}{\pa S_1 \pa S_2}+\frac{\sigma_2^2}{2} \frac{\pa^2 V}{\pa S_2^2} \right)dt\\&
+\left(\sigma_1 \frac{\pa V}{\pa S_1}+\sigma_2 \frac{\pa V}{\pa S_2} \right) dW_t.
\end{split}
\end{equation}

Now, we specify the parameters and principles of the model. 
\subsection{Model parameters}

We use the following notations.

\begin{itemize}

\item $t=0$ is the moment of signing of the coalition agreement and $t=T$ is the moment of elections.

\item $S_i$, $i=1,2$, is level of support to the party $P_i$, $i=1,2$, given in percentage. The variables $S_i$, $i=1,2$, are considered as independent unless we write $S_i(t)$, $i=1,2$. In this case, $S_i(t)$, $i=1,2$, are the stochastic processes which respectively denote the support to the party $P_i$, $i=1,2$, at the moment $t$. In the sequel, $P_1$ denotes the "major" party and $P_2$ denotes the "minor" party.

\item $\mu_i$, $i=1,2$, is an average change in support to the party $P_i$, $i=1,2$, and we consider it constant.

\item $\sigma_i$, $i=1,2$, is the fluctuation of the support to the party $P_i$, $i=1,2$, and we consider it constant. 

\item $E=S_1(0)+S_2(0)$ is an expected level of support to the coalition of parties $P_1$ and $P_2$ at the elections. It is the information obtained from an independent agency at the moment of the conclusion of the coalition agreement. 
 
\item $X$ denotes the election threshold.
 
\item $Y$ denotes value of one mandate i.e. percentage of support necessary to get one mandate based on the poll results at the moment of the coalition contract signing. It is not a constant and it depends on the expected support $E$. We note that larger expected support $E$ implies smaller $Y$ (see the last section). 

\item $Y'=kY-S_2(0)$ for $k$ such that $(k-1)Y-S_2(0)<0$ i.e. it is the percentage of votes that the party $P_2$ needs to get another mandate (if we neglect the election threshold) based on the poll results at the moment of the coalition contract signing.   
 
\item $V(t,S_1,S_2)$  is the "price" of the coalition agreement i.e. the percentage of votes which the "major" party should be willing to transfer to the "minor" party. 

\end{itemize} Regarding the coefficients, we shall need to determine $\sigma_1$, $\sigma_2$ and $T$. The volatilities $\sigma_1$ and $\sigma_2$ can be determined from the historic data while $T$ determines how far in the past we go back. It represents the moment of elections, but it is non-dimensional quantity and it depends on the risk inclination of the "major" party. Strictly rationally speaking, if the "minor" party has a stable support then the "major" party  can allow larger $T$. If the "minor" party does not have a stable support, then $T$ should be set smaller (longer the period is harder is to predict the dynamics of voters' support). However, to mitigate any vagueness, we introduce the principle of equitable treatment for minor parties (refer to point (d) below), stipulating that $T$ should be selected to maximize the proportion of votes ceded by the larger party to the minor one.

Our primary concern lies in determining the value of $V(0, S_1(0), S_2(0))$, as the agreement should be formalized at the initial juncture, typically corresponding to the final day for candidate submissions. In the context of our case study, Montenegro, this moment transpires one month ahead of the elections. If this calculated value, when added to the votes secured by the coalition, meets or surpasses the election threshold needed to secure a single mandate in addition to the projected percentage of votes denoted by $E$, then the "minor" party will be granted an extra mandate. Stated differently, even if the total vote tally falls short of acquiring one more mandate beyond what the projected percentage of votes $E$ offers, yet remains sufficiently proximate to that figure, the "major" party would still concede an additional mandate to the "minor" party.

\begin{example}
\label{ex1}
Let us illustrate this using an example from Table 1. According to pre-electoral polls, Party FfM had a projected support of $S_1(0) = 13.5\%$, while Party JfA had a projected support of $S_2(0) = 1\%$ (although it did not appear in official polls, for simplicity, we consider it as $1\%$). The expected support for the coalition would then be $E = S_1(0) + S_2(0) = 14.5\%$. At this level, approximately $Y=1.15\%$ of the votes correspond to one mandate, and currently, the "minor" party JfA has zero mandates. In the coalition, it would need $Y'=0.15\%$ for one mandate.

Using our model (see \eqref{model}), we calculate $V(0, 13.5, 1)$. If the results on the day of the elections ($t=T$) satisfy the condition $V(0, 13.5, 1) + S_1(T) + S_2(T) \geq E+  0.15\%$, the "minor" party will secure a mandate. Since the actual elections support was $S_1(T) + S_2(T) = 17.53\%$, which is greater than $E + 0.15\% = 14.65\%$, the minor party would obtain the mandate regardless of the value of $V(0, 13.5, 1)$.

A more intriguing situation within our model's perspective involves the relationship between two other ideologically similar parties, which will be discussed in the subsequent sections.
\end{example}

\subsection{Model derivation}
Before we go to the model derivation let us fix several necessary assumptions.

\begin{itemize}

\item[(a)] There exists an independent polling agency accepted by all parties which can provide us with the correct (up to a statistical error) popularity of a party at any time.

\item[(b] Parties have statistically stable support from the moment $t=0$ until the moment of elections in the sense that they will not lose or gain support if no actions are undertaken. This is equivalent to the guaranteed (risk free) interest rate assumption in financing stating that we can always save or loan the money at the guaranteed interest rate. Here, we assume that the support (i.e. the "interest rate") is zero i.e. that we can count on no inflation or deflation of the support.     

\item[(c)] We operate under the premise of a "fair political market," wherein no party can enhance its support without assuming risk. This means that we exclude scenarios akin to financial arbitrage, where risk-free gains can be obtained.

\end{itemize}

Let us now derive the model. We consider the following political portfolio
$$
\Pi=V-\Delta S_1
$$ consisting of one "option" $V$ i.e. the part of the support that the "major" party is willing to transfer and the part $\Delta$ of the support which is spent on $V$. Change of the portfolio in one time step is
$$
d\Pi=dV-\Delta d S_1
$$ and according to \eqref{1} and \eqref{ito}, we have
\begin{equation}
\label{m1}
\begin{split}
d\Pi&=\left( \frac{\pa V}{\pa t}-\Delta \mu_1  + \mu_1 \frac{\pa V}{\pa S_1}+\mu_2 \frac{\pa V}{\pa S_2}+\frac{\sigma_1^2}{2} \frac{\pa^2 V}{\pa S_1^2}+{\sigma_1\sigma_2} \frac{\pa^2 V}{\pa S_1 \pa S_2}+\frac{\sigma_2^2}{2} \frac{\pa^2 V}{\pa S_2^2} \right)dt\\&
+\left(\sigma_1 \frac{\pa V}{\pa S_1}+\sigma_2 \frac{\pa V}{\pa S_2}-\Delta \sigma_1  \right) dW_t.
\end{split}
\end{equation} Now, we want to eliminate the random component in the latter equation and we take
\begin{equation*}
\begin{split}
&\sigma_1 \frac{\pa V}{\pa S_1}+\sigma_2 \frac{\pa V}{\pa S_2}-\Delta \sigma_1=0 \ \ \implies \\
&\Delta = \frac{\pa V}{\pa S_1}+\frac{\sigma_2}{\sigma_1} \frac{\pa V}{\pa S_2}.
\end{split}
\end{equation*} Thus, \eqref{m1} becomes
\begin{equation}
\label{m2}
\begin{split}
d\Pi&=\left( \frac{\pa V}{\pa t} +\big(\mu_2-\mu_1 \frac{\sigma_2}{\sigma_1} \big) \frac{\pa V}{\pa S_2}+\frac{\sigma_1^2}{2} \frac{\pa^2 V}{\pa S_1^2}+{\sigma_1\sigma_2} \frac{\pa^2 V}{\pa S_1 \pa S_2}+\frac{\sigma_2^2}{2} \frac{\pa^2 V}{\pa S_2^2} \right)dt.
\end{split}
\end{equation} Next, we denote by $X$ the election threshold and consider the latter equation in the following domain
\begin{equation}
\label{dom}
\begin{split}
(S_1,S_2)\in D:=(X,40)\times (0,X) 
\end{split}
\end{equation} i.e. $S_1$ will always lie above the election threshold and $S_2$ will always be under the election threshold. We also assume that $S_1$ is constantly under $40\%$ since otherwise, there is no sense to make a coalition. We note that in a system with lots of parties in which case the current research has sense, a party with over $40\%$ safely wins the elections and it needs no coalitions. Under such assumption, we can take that 
$$
\frac{\pa V}{\pa S_2} \sim 0
$$ i.e. the value $V$ will not change significantly with the change of $ S_2 $ since the change of $S_2$ itself is small. Thus, \eqref{m2} becomes 
\begin{equation}
\label{m3}
\begin{split}
d\Pi&=\left( \frac{\pa V}{\pa t} + \frac{\sigma_1^2}{2} \frac{\pa^2 V}{\pa S_1^2}+{\sigma_1\sigma_2} \frac{\pa^2 V}{\pa S_1 \pa S_2}+\frac{\sigma_2^2}{2} \frac{\pa^2 V}{\pa S_2^2} \right)dt.
\end{split}
\end{equation} Next, taking into account assumptions (b) and (c) from the above, we see that 
$$
d\Pi=0
$$ since the party cannot lose support according to (b), but it cannot also gain it according to (c) (i.e. change in the political portfolio must be zero in the arbitrage free world). From here and \eqref{m3}, we finally conclude 
\begin{equation}
\label{model}
\begin{split}
0&= \frac{\pa V}{\pa t} + \frac{\sigma_1^2}{2} \frac{\pa^2 V}{\pa S_1^2}+{\sigma_1\sigma_2} \frac{\pa^2 V}{\pa S_1 \pa S_2}+\frac{\sigma_2^2}{2} \frac{\pa^2 V}{\pa S_2^2} .
\end{split}
\end{equation}This is a backward linear degenerate parabolic equation and we need to impose initial and boundary conditions in order to get existence of a unique solution $V$. We set (recall that $Y'$ is the percentage of votes that $P_2$ needs for another mandate based on the polling at the moment of the coalition contract signing)
\begin{equation}
\label{bc-ic}
\begin{split}
&V(T,S_1,S_2)=\max\{S_1+S_2-E-Y',0 \} \ \ \text{(final condition/moment of elections);}\\
&V(t,S_1,0)=V(t,S_1,X)=V(t,X,S_2)=V(t,40,S_2)=0 \ \ \text{(boundary conditions)}.
\end{split}
\end{equation} Let us explain the latter conditions.

The final condition states that if the result of the elections is larger than the expected one for the value of the missing support $Y'$, than the "minor" party should get another mandate. Since we tacitly assume $S_1+S_2 \sim S_1$ due to difference in the strength of the parties, we say that the "major" party renounced $Y'$ percent of votes in favor of the minor party. If the result is smaller than $E+Y'$ then the "major" party has no spare votes. 

The condition $V(t,S_1,0)=0$ means that if strength of the "minor" party is zero, then there is no sense for the "major" party to enter the coalition and give away its votes. 
 
If the "minor" party reached the election threshold, then there is no need for additional votes from the "major" party and we set $V(t,S_1,X)=0$.

If the "major" party has merely reached the election threshold, then it has no additional votes to give away. Therefore, $V(t,X,S_2)=0$.

Finally, if the "major" party reached absolute support which realistically lies on $40\%$ threshold, then there is no motivation for entering coalitions and we have $V(t,40,S_2)=0$.

The meaning of the model is that the "major" party would surely agree with one more mandate for the "minor" party if the elections results are for $Y'$ greater than the expected results. 


Now, we must address the question of determining the appropriate value for the final time, denoted as $T$. The coalition agreement is finalized one month prior to the elections. Nevertheless, converting this one-month period into the final time $T$ poses a challenge since we are in a non-dimensional setting. In this pursuit, we delve into the characteristics of the equation itself. This equation is a variation of the backward heat equation. The final conditions mirror the heat distribution at the moment $T$, while the zero boundary conditions signify a continual cooling of the space back in time.

In simpler terms, when considering an exceedingly large value for $T$, the solution's value $V$ at $t=0$ also converges to zero. This is because the cooling at the boundaries eliminates the impact of the initial heat distribution. Consequently, we propose the following approach:

\begin{itemize}

\item[(d)] Principle of Fair Treatment for the "Minor" Party. We shall select $T$ based on the condition that maximizes the value $V(0, S_1(0), S_2(0))$, where $S_1(0)$ and $S_2(0)$ denote the predicted support levels of the two parties at the moment the coalition agreement is signed. This principle is aimed at ensuring equitable treatment for the party with lesser influence.

\end{itemize}

In other words, we consider problem \eqref{model} for $t\leq 0$ and take the following final-boundary conditions
\begin{equation}
\label{bc-ic-true}
\begin{split}
&V(0,S_1,S_2)=\max\{S_1+S_2-E-Y',0 \} \ \ \text{(final condition/moment of elections);}\\
&V(t,S_1,0)=V(t,S_1,X)=V(t,X,S_2)=V(t,40,S_2)=0 \ \ \text{(boundary conditions)}.
\end{split}
\end{equation} Then, we find $T_{max}>0$ such that 
\begin{equation}
\label{max}
V(-T_{max},S_1(0),S_2(0)) = \max\limits_{t\geq 0} V(-t,S_1(0),S_2(0)).
\end{equation} The latter quantity will represent the percentage of votes that the "major" party is willing to give away to the "minor" one. Solving \eqref{model}, \eqref{bc-ic-true} and taking $V(-T_{max},S_1(0),S_2(0))$ for the final result is equivalent as solving \eqref{model}, \eqref{bc-ic} with $T=T_{max}$ and taking $V(0,S_1(0),S_2(0))$ as the final result.

\subsection{Solving the problem \texorpdfstring{\eqref{model}, \eqref{bc-ic}}{Lg}}

We note that equation \eqref{model} can be rewritten in the form
$$
\frac{\pa V}{\pa t}+\frac{1}{2} \left(\sigma_1 \frac{\pa}{\pa S_1}+\sigma_2 \frac{\pa}{\pa S_2} \right)^2V=0.
$$ Thus, introducing the change of variables
$$
y_1=S_1, \ \ y_2=\sigma_1 S_1+\sigma_2 S_2
$$ equation \eqref{model} becomes (we keep $V$ for the unknown function depending now on $y_1$ and $y_2$)
\begin{equation}
\label{model-1}
\pa_tV+\frac{1}{2} \pa_{y_2 y_2} V=0,
\end{equation} while the domain $D$ in which we solve the equation \eqref{model} is transformed into
$$
D'=\{(y_1,y_2):\, X\leq y_1 \leq 40, \, \sigma_1 y_1 \leq y_2 \leq \sigma_1 y_1+\sigma_2 X \}.
$$ The final data remain the same and the boundary data are zero. Since \eqref{model-1} does not depend on $y_1$, we can consider it as a parameter. The standard separation of variables gives an explicit solution in the form
\begin{equation}
\label{expl-sol}
V(t,y_1,y_2)=\sum\limits_{n=1}^\infty A_n(y_1)\sin\left(\frac{n\pi}{\sigma_2 X} (y_2-{\sigma_1 y_1})\right) e^{-\frac{1}{2}\left(\frac{n \pi}{\sigma_2 X}\right)^2 (T-t)}
\end{equation} where
$$
A_n(y_1)=\frac{1}{\sigma_2 X} \int_{\sigma_1 y_1}^{\sigma_1 y_1+\sigma_2 X} \sin\left(\frac{n\pi}{\sigma_2 X} y_2 \right) V(T,y_1,y_2)dy_2.
$$ However, the problem is easily solved using the finite difference method which will be demonstrated on the concrete examples in the sequel.

\section{Elections in Montenegro -- a case study}

We shall now analyze a few examples from practice. The parliament in Montenegro has 81 representative and the election threshold is $3\%$. Recent elections in Montenegro are convenient since fifteen parties and alliances took part in the elections and two of them were very close to the election threshold.  If they made a pre-elections agreement with ideologically similar parties, they would have the parliamentary status and the partner parties would even have more representatives (from purely theoretical point of view, i.e., by adding supports obtained by the respective parties and disregarding potential loss/gain of votes due to entering the coalition).

We recall that pre-elections data are taken from \cite{CEDEM} and the elections results from \cite{DIK}.

Let us first very briefly introduce the different political blocks in Montenegro. While the following description may not be entirely precise, its accuracy is inconsequential to our theoretical exploration of election strategies. We also stress again that, in considering potential coalition, we have simply sum up votes that the parties obtained as individual entities taking part in the elections.
 \begin{itemize}
 \item Right-wing pro-serbian parties and coalitions characterized by, roughly speaking, the following: support Russia, do not recognize Kosovo independence, do not recognize the genocide in Srebrenica, and they are anti-Nato. We mention only those of our interest: For the Future of Montenegro (FfM), Socialist People Party (SPP), and Justice For All (JfA).
 
 \item Right-wing Montenegrin parties and coalitions characterized by the narrative that Montenegro is directly endangered by Serbia, strongly support NATO and closely follow the EU foreign policy. We mention here only Democratic Party of Socialists (DPS) and Social Democratic Party (SDP).
 
\item Neutral parties supporting NATO and EU whose narrative is mainly focused on economy and fight against crime. We mention Movement Europe Now (MEN), Democrats (DEM), and United Reform Action (URA).

\item Minority parties of Bosnian (BS), Albanians (Al), and Croats (Cro) with politics close to the neutral parties.  
 \end{itemize} Interestingly, although we have fairly polarized political situation in Montenegro, unlike theoretical predictions \cite{golder}, there is still an abundance of various political entities taking part in the elections process without entering larger coalitions.

 In Example \ref{ex1}, we explained the pre and post electoral support of pro-serbian parties FfM and JfA. We have seen that pre-electoral coalition would be of great interest for both FfM and JfA since they would both get one more representative and MEN and DPS would lose one representative each (see Table 1). Here, we shall compute $V(0,13.5,1)$ (where $13.5\%$ and $1\%$ are pre-electoral support of FfM and JfA). Knowing $V(0,13.5,1)$, the parties should make the following agreement: {\em if the total elections support $S_1+S_2$ of the coalition exceeds $E+Y'=14.5+0.15=14.65$, then JfA is entitled to one mandate}. 
 
\begin{figure}[h!]
  \centering
\stackunder[3pt]{\includegraphics[width=6in,height=2in]{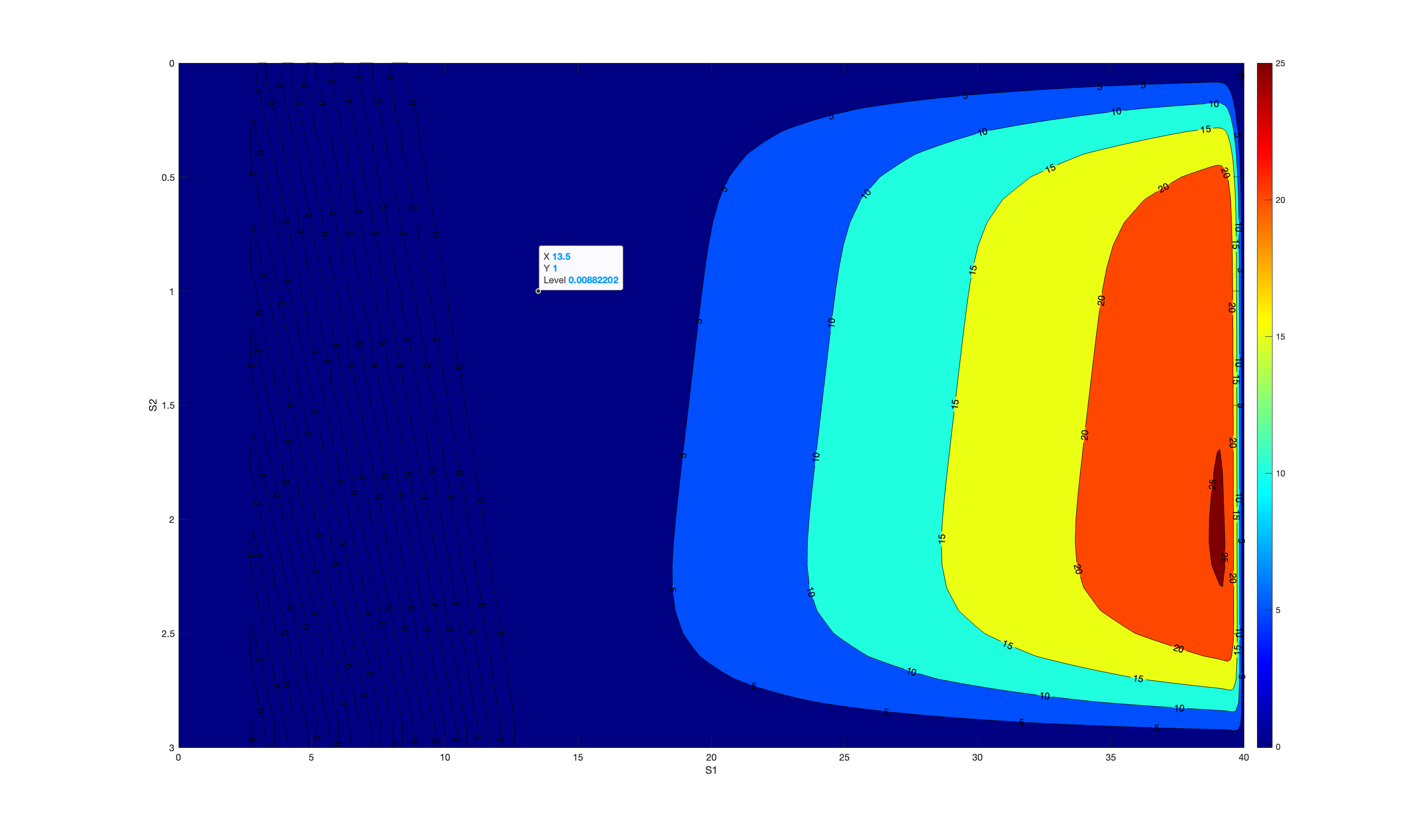}}{$V(-100,13.5,1)=0.0088$}  

\stackunder[3pt]{\includegraphics[width=6in,height=2in]{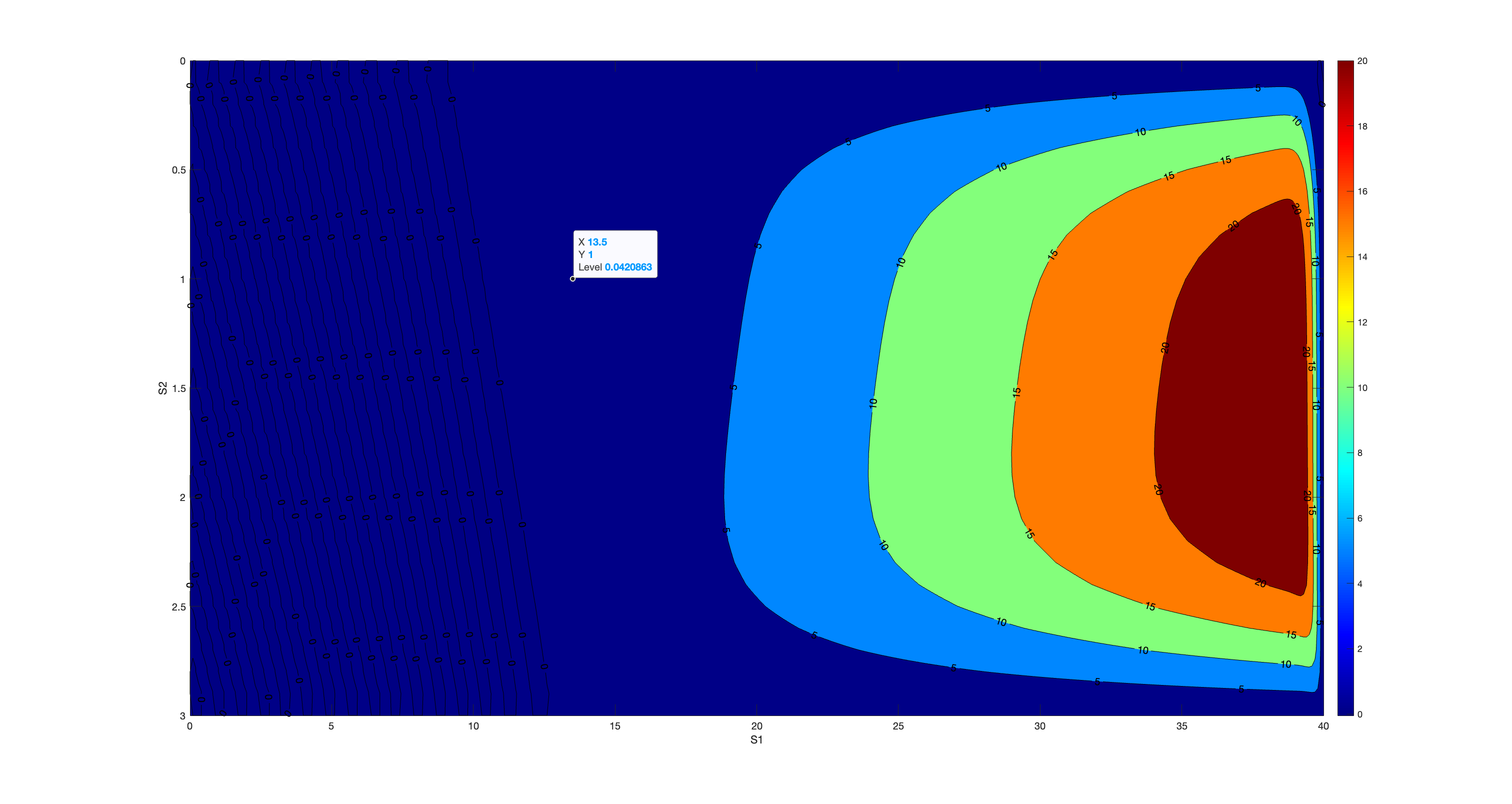}}{$V(-200,13.5,1)=0.042$}  
  
\stackunder[3pt]{\includegraphics[width=6in,height=2in]{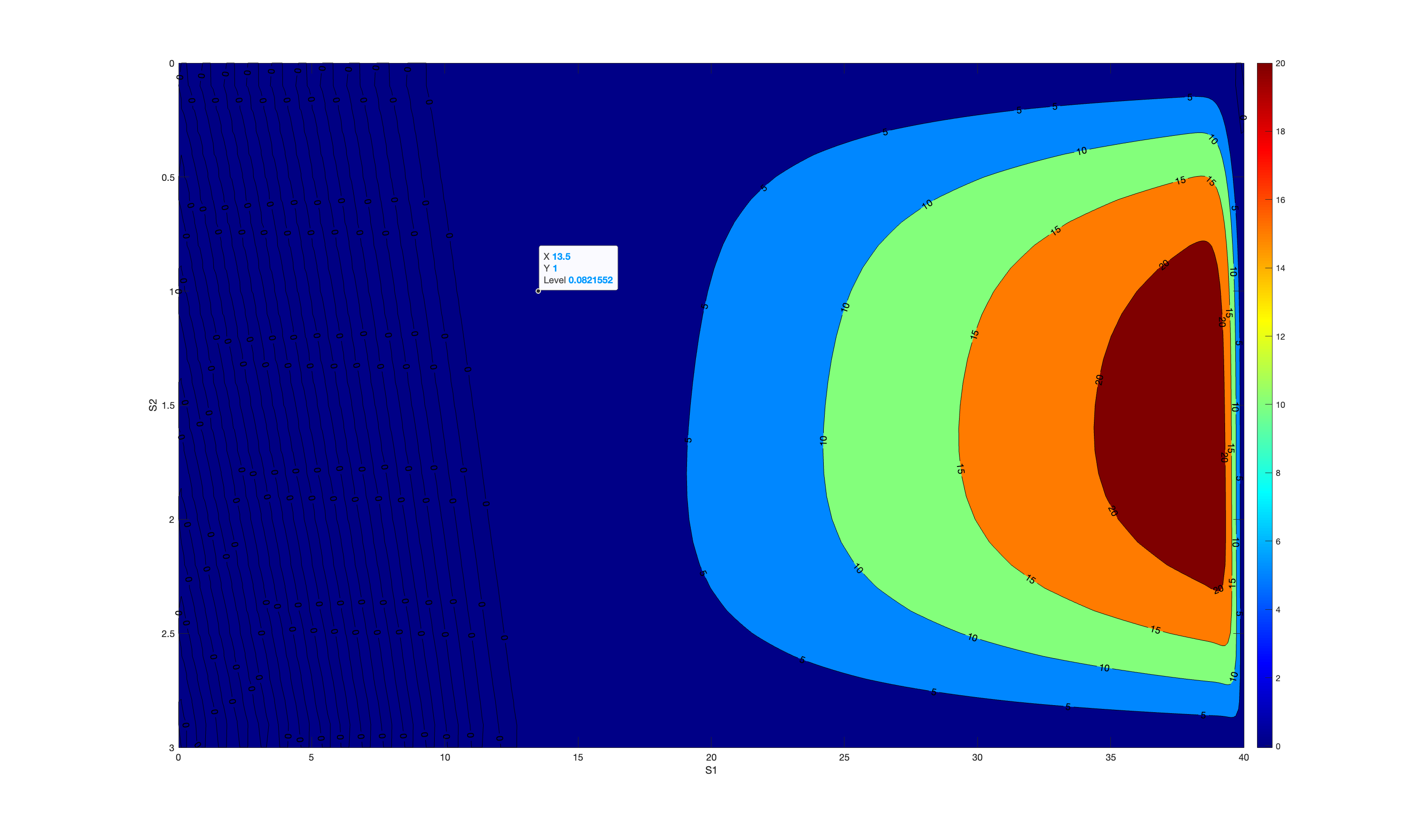}}{$V(-300,13.5,1)=0.082$}  

\stackunder[3pt]{\includegraphics[width=6in,height=2in]{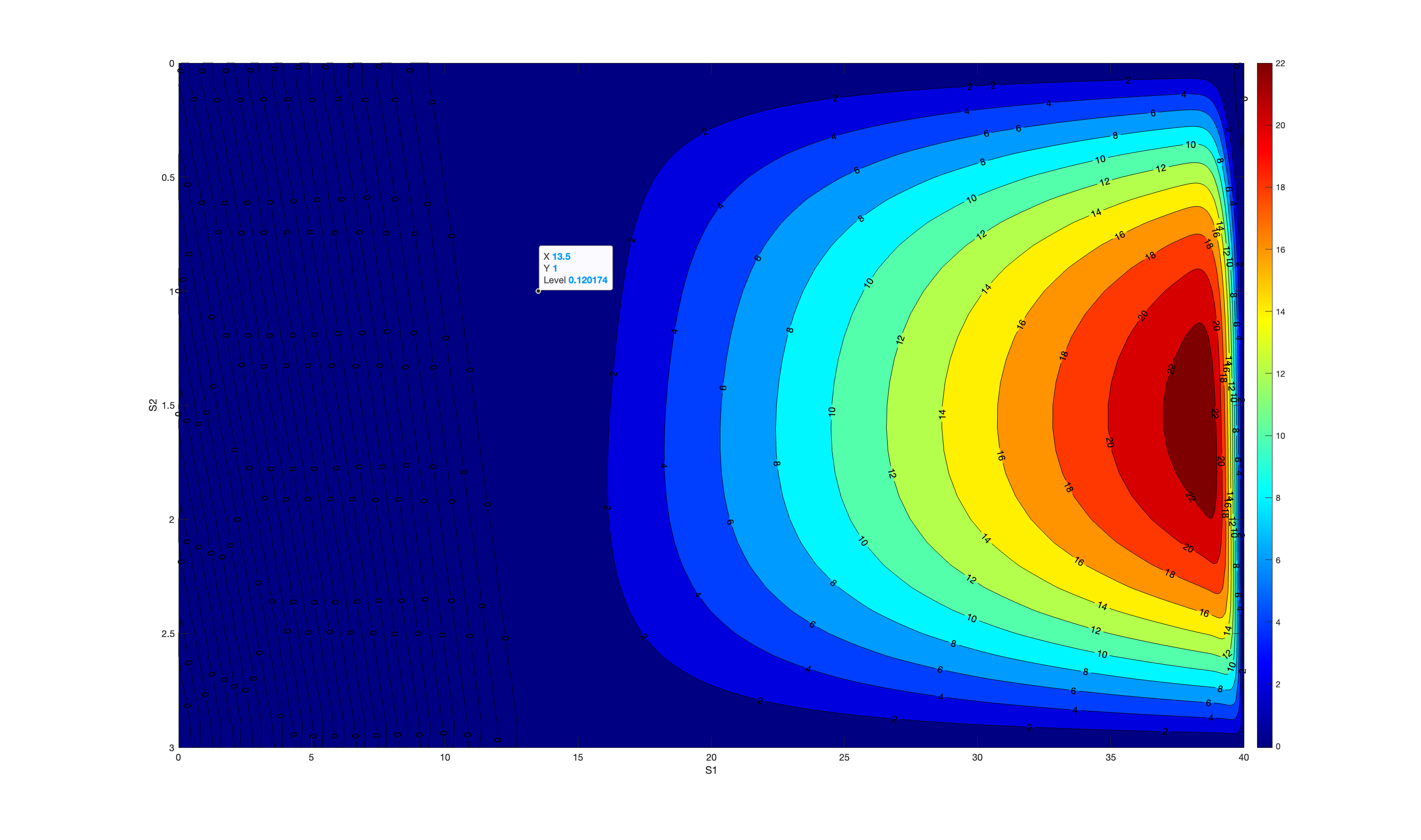}}{$V(-400,13.5,1)=0.1201$}  

\caption{Values for $T=100, 200, 300, 400$ for a potential coalition FfM and JfA.}

  \label{figFfM}
\end{figure}

 Although in this concrete situation we have $S_1+S_2=17.53>E+Y'$ i.e. JfA would get a mandate independently of the value of $V(0,13.5,1)$, let us still show how to compute the latter. We do not have concrete historical data and we shall take $\sigma_1=\sigma_2=1$ i.e. volatility is fairly small (it seems that it is often a case for conservative parties on Balkan although we do not know actual data). For the final data, we shall use the principle of fair treatment of the "minor" party given by (d) above. With the latter assumptions at hand, we actually aim to solve the following special case of \eqref{model}, \eqref{bc-ic}
\begin{equation}
\label{model-2}
\begin{split}
0&= \frac{\pa V}{\pa t} + \frac{1}{2} \frac{\pa^2 V}{\pa S_1^2}+ \frac{\pa^2 V}{\pa S_1 \pa S_2}+\frac{1}{2} \frac{\pa^2 V}{\pa S_2^2} \\
&V(0,S_1,S_2)=\max\{S_1+S_2-15.65,0 \} \ \ \text{(final condition/moment of elections);}\\
&V(t,S_1,0)=V(t,S_1,3)=V(t,3,S_2)=V(t,40,S_2)=0 \ \ \text{(boundary conditions)}.
\end{split}
\end{equation} The results are given in Figure \ref{figFfM} for $T=100, 200, 300, 400$ iterations (i.e. assuming that the elections happened after $100, 200, 300$ or $400$ iterations). We did not look for the exact $T_{max}$ (given by \eqref{max}) since $S_1(T)+S_2(T)+V(-T,S_1(0),S_2(0))>S_1(T)+S_2(T)>15.65$ for any chosen $T>0$ representing the moment of elections. Thus, JfA gets one mandate.

\begin{figure}[h!]
  \centering
\stackunder[3pt]{\includegraphics[width=6in,height=2in]{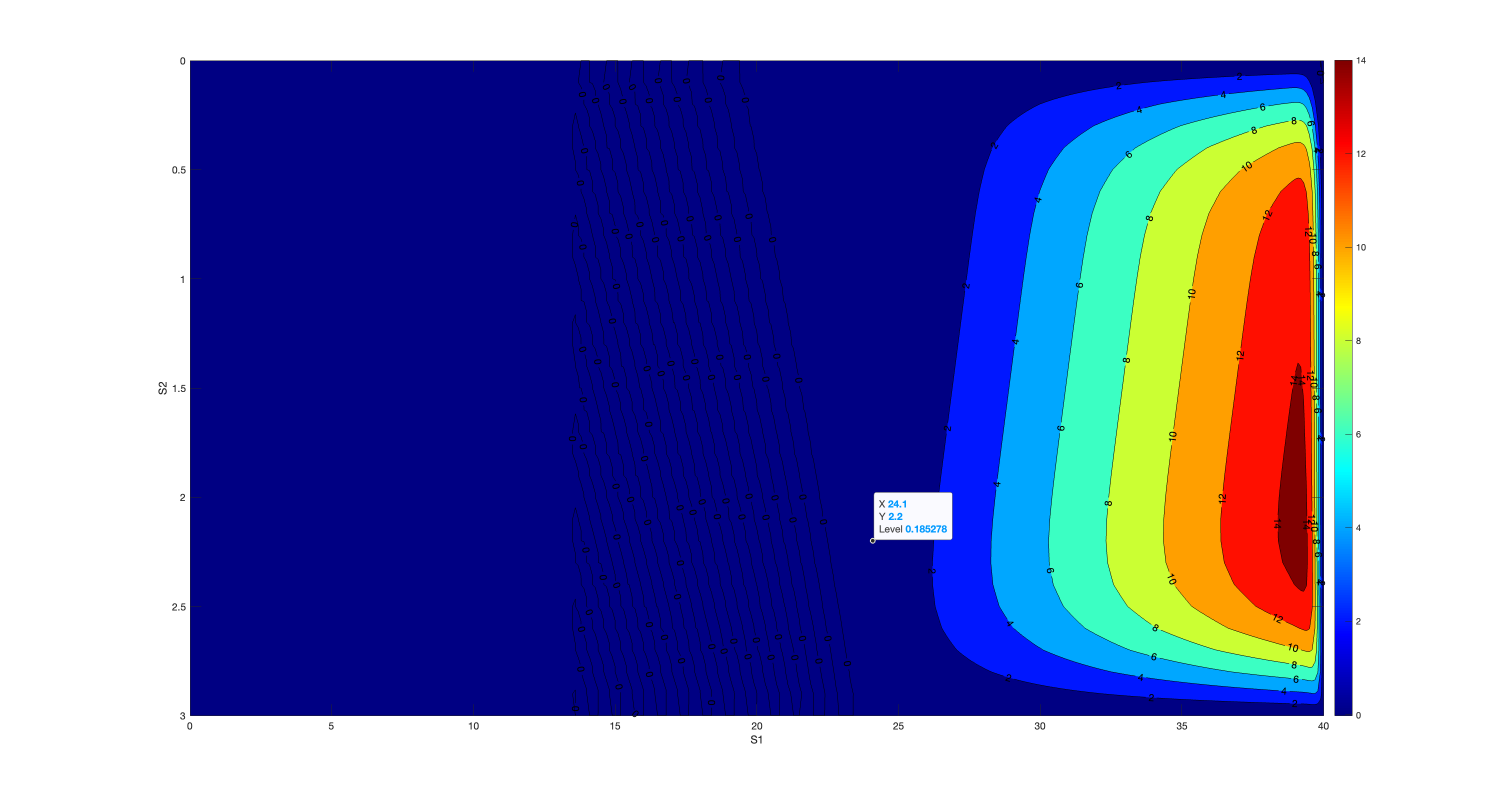}}{$V(-100,24.1,2.2)=0.1853$}  

\stackunder[3pt]{\includegraphics[width=6in,height=2in]{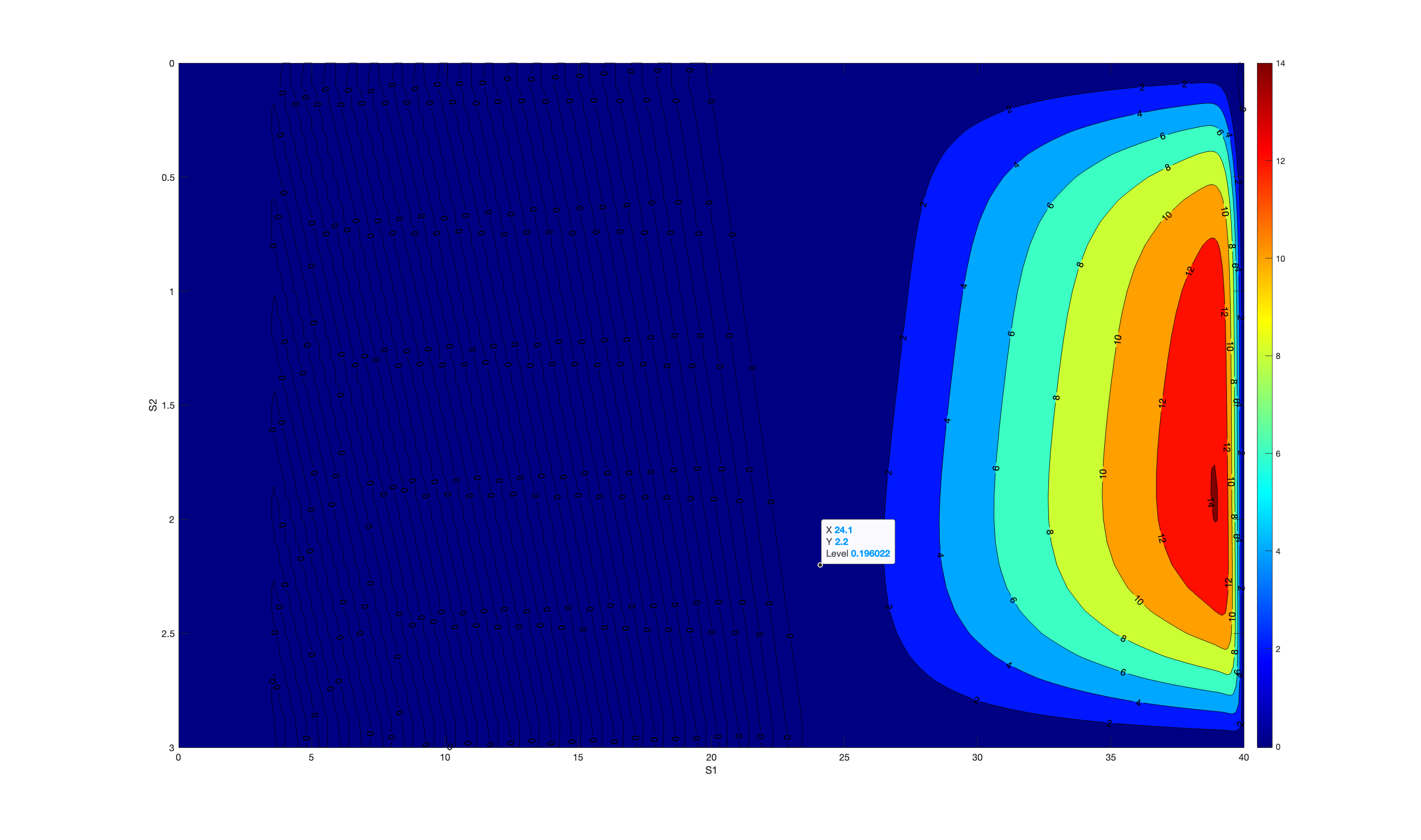}}{$V(-200,24.1,2.2)=0.196$}  
  
\stackunder[3pt]{\includegraphics[width=6in,height=2in]{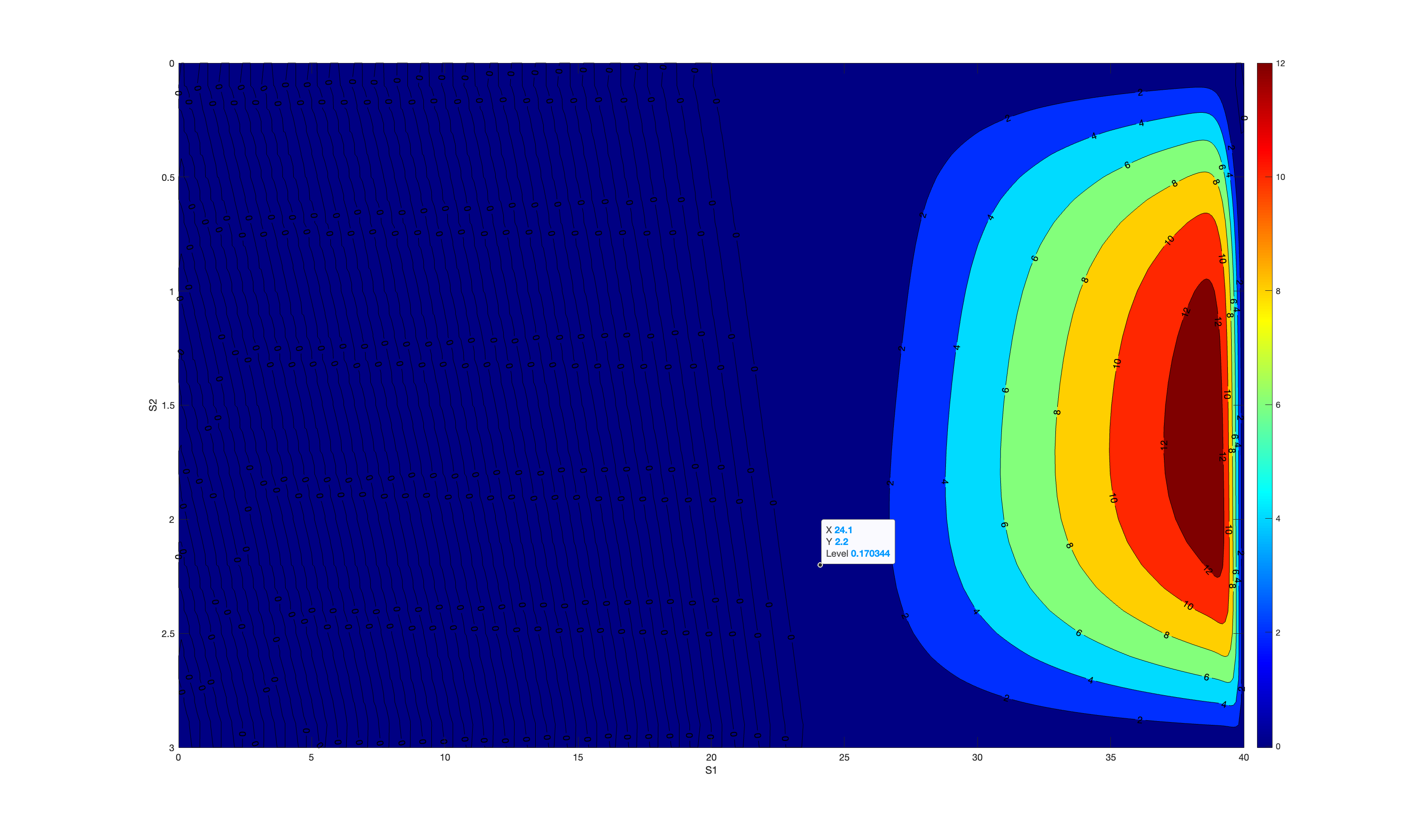}}{$V(-300,24.1,2.2)=0.1703$}  

\stackunder[3pt]{\includegraphics[width=6in,height=2in]{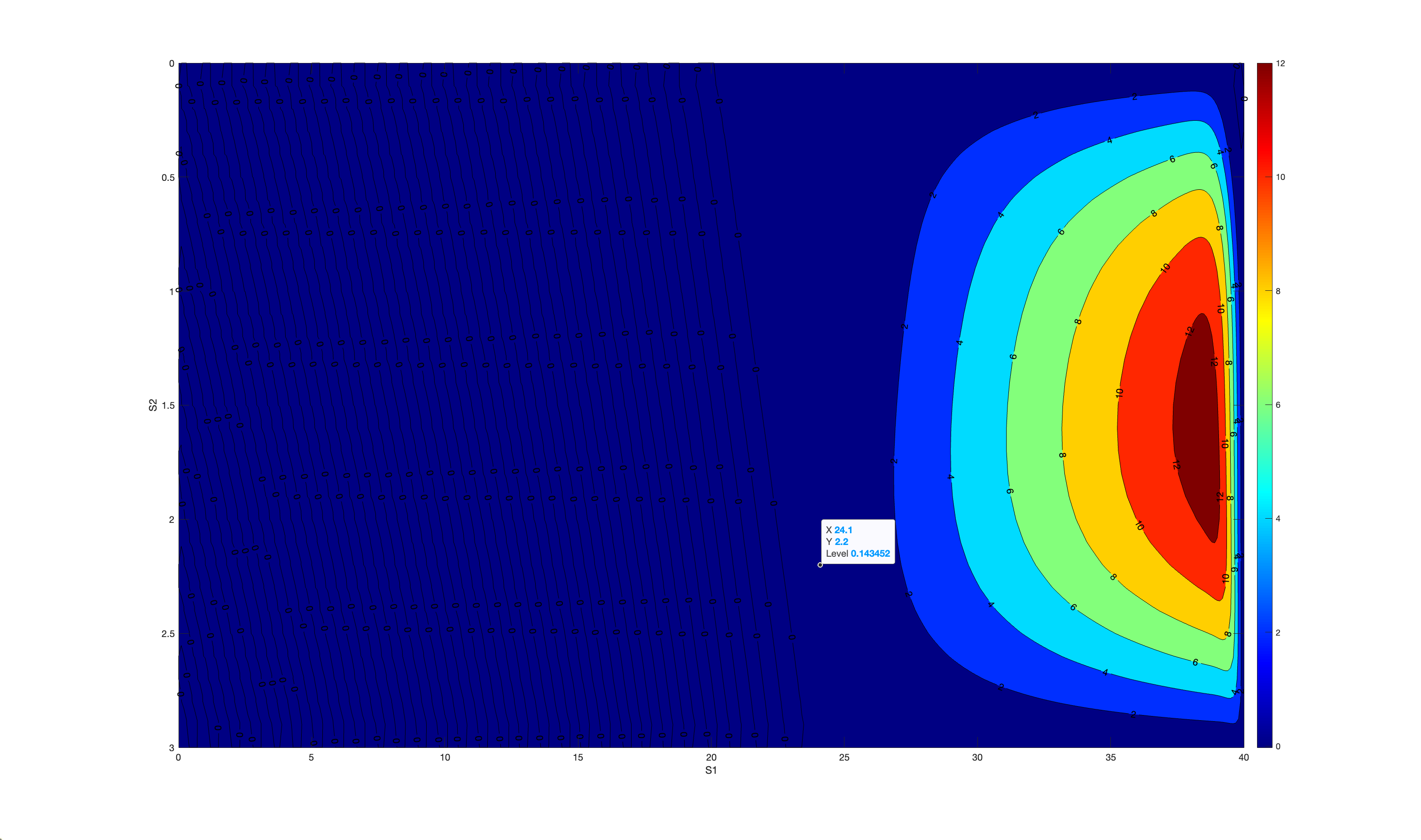}}{$V(-400,24.1,2.2)=0.1434$}  

\caption{Values for $T=100, 200, 300, 400$ for a potential coalition DPS and SDP.}

  \label{figDPS}
\end{figure}

The second example from the same elections is more interesting since it is not clear whether the "minor" party would get another mandate without computing $V(0,S_1(0),S_2(0))$. To be more precise, let us consider the results of DPS and SDP (Montenegrin right-oriented parties). Pre-elections polls gave DPS the result $S_1(0)=24.1\%$ and SDP the result $S_2(0)=2.2\%$. In the coalition, for one mandate one would need approximately $1.15\%$. This means that SDP misses only $0.1\%$ for two mandates. Actual elections results were $S_1(T_{max})+S_2(T_{max})=26.21\%$.
Thus, if 
$$
S_1(T_{max})+S_2(T_{max})+V(-T_{max},24.1,2.2)>E+Y'=26.4\%
$$then SDP gets two mandates. Otherwise it has one mandate. To find out this, we again take $\sigma_1=\sigma_2=1$, and solve

\begin{equation}
\label{model-3}
\begin{split}
0&= \frac{\pa V}{\pa t} + \frac{1}{2} \frac{\pa^2 V}{\pa S_1^2}+ \frac{\pa^2 V}{\pa S_1 \pa S_2}+\frac{1}{2} \frac{\pa^2 V}{\pa S_2^2} \\
&V(0,S_1,S_2)=\max\{S_1+S_2-26.4,0 \} \ \ \text{(final condition/moment of elections);}\\
&V(t,S_1,0)=V(t,S_1,3)=V(t,3,S_2)=V(t,40,S_2)=0 \ \ \text{(boundary conditions)}.
\end{split}
\end{equation} The results are provided in Figure \ref{figDPS} for $T=100, 200, 300, 400$ iterations. We see that the maximum is reached for $T=200$ i.e. $V(-200,S_1(0),S_2(0))=0.196$. Thus $S_1(200)+S_2(200)+V(-200,24.1,2.2)=26.41>26.4$ and SDP gets two mandates.

\section{Data availability}

Raw data are available at the web pages of the Montenegro State Election Committee: \url{dik.co.me}, and non-governmental organisation CEDEM: \url{cedem.me}. Derived data supporting the findings of this study are available from the corresponding author of this paper on request.

\section{Acknowledgment}
This research was funded in full by the Austrian Science Fund (FWF) under the grant number {\bf P 35508}. For the purpose of open access, the author has applied a CC BY public copyright licence to any Author Accepted Manuscript version arising from this submission.

I would like to thank my friend and colleague prof. Andrej Novak from the University of Zagreb for working out the numerical simulations in the paper.


\end{document}